\documentclass[amsfonts, 12pt]{article}

\usepackage{amsmath,amsthm,amsfonts,amssymb}

\newtheorem{theorem}{Theorem}[section]

\newtheorem{lemma}[theorem]{Lemma}

\newtheorem{remark}[theorem]{Remark}

\def\cA{\mathcal{A}}
\def\cB{\mathcal{B}}

\def\cH{\mathcal{H}}
\def\cF{\mathcal{F}}

\def\cP{\mathcal{P}}
\def\cS{\mathcal{S}}

\def\bC{\mathbb{C}}

\def\bE{\mathbb{E}}
\def\bR{\mathbb{R}}

\begin{document}

\title{It\^o formula for integral processes related to space-time L\'evy white noise}

\author{Raluca M. Balan\footnote{Corresponding author. University of Ottawa, Department of Mathematics and Statistics,
585 King Edward Avenue, Ottawa, ON, K1N 6N5, Canada. E-mail
address: rbalan@uottawa.ca} \footnote{Research supported by a
grant from the Natural Sciences and Engineering Research Council
of Canada.}
\and
Cheikh B. Ndongo\footnote{University of Ottawa, Department of Mathematics and Statistics,
585 King Edward Avenue, Ottawa, ON, K1N 6N5, Canada. E-mail
address: cndon072@uottawa.ca}}

\date{May 18, 2015}
\maketitle

\begin{abstract}
\noindent In this article, we give a new proof of the It\^o formula for some integral processes related to the space-time L\'evy noise introduced in \cite{B15} and \cite{B15-stable} as an alternative for the Gaussian white noise perturbing an SPDE.
We discuss two applications of this result, which are useful in the study of SPDEs driven by a space-time L\'evy noise with finite variance: a maximal inequality for the $p$-th moment of the stochastic integral, and the It\^o representation theorem leading to a chaos expansion similar to the Gaussian case.
\end{abstract}

\vspace{10mm}


\vspace{5mm}

%60H05 stochastic integrals
%60G51 processes with independent increments and Levy processes

\pagebreak

\section{Introduction}

Random processes indexed by sets in the space-time domain are useful objects in stochastic analysis, since they can be viewed as mathematical models for the noise perturbing a stochastic partial differential equation (SPDE). In the recent years, a lot of effort has been dedicated to studying the behaviour of the solution of basic equations (like the heat or wave equations), driven by a
{\em Gaussian white noise}. This type of noise was introduced by Walsh in \cite{walsh86} and is defined as a zero-mean Gaussian process $W=\{W(B); B \in \cB_b(\bR_{+} \times \bR^d)\}$,  with covariance $E[W(A)W(B)]=|A \cap B|$, where $|\cdot|$ denotes the Lebesgue measure and
$\cB_b(\bR_{+} \times \bR^d)$ is the class of bounded Borel sets in $\bR_{+} \times \bR^d$.

In the recent articles \cite{B15} and \cite{B15-stable}, a new process has been introduced as an alternative for the Gaussian white noise perturbing an SPDE, which has a structure similar to a L\'evy process. We introduce briefly the definition of this process below.

Let $N$ be a Poisson random measure (PRM) on $\bE=\bR_{+} \times \bR^d \times \bR_0$
of intensity $\mu=dtdx \nu(dz)$ where $\bR_0=\bR \verb2\2\{0\}$ and $\nu$ is a L\'evy measure on $\bR$:
$$\int_{\bR_0}(1 \wedge |z|^2)\nu(dz)<\infty \quad \mbox{and} \quad \nu(\{0\})=0.$$
We denote by $\widehat{N}$ the compensated PRM defined by $\widehat{N}(A)=N(A)-\mu(A)$ for any Borel set $A$ in $\bE$ with $\mu(A)<\infty$.
The L\'evy-type noise process mentioned above is defined as $Z=\{Z(B); B \in \cB_b(\bR_{+} \times \bR^d)\}$, where
$$Z(B)=a|B|+\int_{B \times \{|z|>1\}}z N(ds,dx,dz)+\int_{B \times \{|z| \leq 1\}}z \widehat{N}(ds,dx,dz),$$
for some $a \in \bR$.
It was shown in \cite{B15-stable} that $Z$ is an ``independently scattered random measure'' (in the sense of \cite{RR89}) with characteristic function:
$$E(e^{iu Z(B)})=\exp\left\{|B|\left(a+\int_{\bR_0}(e^{iuz}-1-iuz1_{\{|z| \leq 1\}})\nu(dz)\right)\right\}, \quad u \in \bR.$$
(In particular, $Z$ can be an $\alpha$-stable random measure  with $\alpha \in (0,2)$, as in Definition 3.3.1 of \cite{ST94}.)
One can define the stochastic integral of a process $X=\{X(t,x);t \geq 0,x \in \bR^d\}$ with respect to $Z$ and for a {\em certain}  integrands, %this integral satisfies:
\begin{eqnarray*}
\lefteqn{\int_0^T \int_{\bR^d}X(t,x)Z(dt,dx)=a\int_0^T \int_{\bR^d}X(t,x) dtdx+}\\
& & \int_0^T\int_{\bR^d}\int_{\{|z|>1\}}X(t,x)zN(dt,dx,dz)+\int_{0}^T \int_{\bR^d}\int_{\{|z|\leq 1\}}X(t,x)z\widehat{N}(dt,dx,dz).
\end{eqnarray*}

The stochastic integral with respect to $\widehat{N}$ (or $N$) can be defined using classical methods (see e.g. \cite{applebaum09}). We review briefly this definition here.

Assume that $N$ is defined on a probability space $(\Omega, \cF,P)$. On this space, we consider the filtration
$$\cF_t= \sigma(\{N([0,s] \times B \times \Gamma); 0 \leq s \leq t, B \in \cB_b(\bR^d), \Gamma \in \cB_b(\bR_0)\}),$$
where $\cB_b(\bR^d)$ is the class of bounded Borel sets in $\bR^d$ and $\cB_b(\bR_0)$ is the class of Borel sets in $\bR_0$ which are bounded away from $0$.

 An {\em elementary} process on $\Omega \times \bR^d \times \bR_0$ is a process of the form
$$H(\omega,t,x,z)=X(\omega)1_{(a,b]}(t)1_{A}(x) 1_{\Gamma}(z),$$
where $0 \leq a<b$, $X$ is an $\cF_a$-measurable bounded random variable, $A \in \cB_b(\bR^d)$ and $\Gamma \in \cB_b(\bR_0)$.
A process $H=\{H(t,x,z); t \geq 0, x \in \bR^d,z \in \bR_0\}$ is called {\em predictable} if it is measurable with respect to the $\sigma$-field $\cP_{\Omega \times \bR_{+} \times \bR^d \times \bR_0}$ generated by all linear combinations of elementary processes.

As in It\^o's classical theory, for any predictable process $H$ such that
\begin{equation}
\label{cond-H-R0}
E\int_0^t \int_{\bR^d} \int_{\bR_0} |H(s,x,z)|^2 \nu(dz)dxds<\infty \quad
\mbox{for all} \ t>0,
\end{equation}
we can define the stochastic integral of $H$ with respect to $\widehat{N}$ and
the process
$\{\int_0^t \int_{\bR^d}\int_{\bR_0}H(s,x,z)\widehat{N}(ds,dx,dz);t \geq 0\}$
is a zero-mean square-integrable martingale which satisfies
\begin{equation}
\label{E-int-hatN}
E\left|\int_0^t \int_{\bR^d}\int_{\bR_0}H(s,x,z)\widehat{N}(ds,dx,dz)\right|^2=\int_0^t \int_{\bR^d}\int_{\bR_0}|H(s,x,z)|^2 \nu(dz)dxds.
\end{equation}

On the other hand, for any predictable process $K$ such that
$$E\int_{0}^t \int_{\bR^d}\int_{\bR_0}|K(s,x,z)|\nu(dz)dxds<\infty \quad \mbox{for all} \quad t>0,$$
we can define the integral of $K$ with respect to $N$ and this integral satisfies
\begin{equation}
\label{E-int-N}
E\int_0^t \int_{\bR^d}\int_{\bR_0}K(s,x,z)N(ds,dx,dz)=E\int_{0}^t \int_{\bR^d}\int_{\bR_0}K(s,x,z)\nu(dz)dxds.
\end{equation}

In this article, we work with processes whose trajectories are c\`adl\`ag, i.e.
right-continuous with left limits. If $x$ is a c\`adl\`ag function, we denote by $x(t-)=\lim_{s \uparrow t}x(s)$ the left limit at time $t$ and $\Delta x(t)=x(t)-x(t-)$ the jump size at time $t$.
We will prove the following result.

%The following result is the analogue of Theorem 4.4.7 of \cite{applebaum09}.
\begin{theorem}[It\^o Formula I]
\label{Ito-formula-th}
Let $Y=\{Y(t)\}_{t \geq 0}$ be a process defined by
\begin{eqnarray}
\label{integral-proc}
Y(t)&=&\int_0^t G(s)ds+\int_0^t \int_{\bR^d}\int_{\{|z|>1\}}K(s,x,z)N(ds,dx,dz)+\\
\nonumber
& & \int_0^t \int_{\bR^d}\int_{\{|z| \leq 1\}}H(s,x,z)\widehat{N}(ds,dx,dz), \quad t \geq 0,
\end{eqnarray}
where $G$,$K$ and $H$ are predictable processes which satisfy
\begin{eqnarray}
\label{cond-G}
& & E\int_0^t|G(s)|ds<\infty \quad \mbox{for all} \quad t>0, \\
\label{cond-K}
& & E\int_0^t \int_{\bR^d} \int_{\{|z>1\}} |K(s,x,z)| \nu(dz)dxds<\infty \quad
\mbox{for all} \ t>0,\\
\label{cond-H}
& & E\int_0^t \int_{\bR^d} \int_{\{|z|\leq 1\}} |H(s,x,z)|^2 \nu(dz)dxds<\infty \quad
\mbox{for all} \ t>0.
\end{eqnarray}
Then there exists a c\`adl\`ag modification of $Y$ (denoted also by $Y$) such that for any function $f \in C^2(\bR)$ and for any $t>0$, with probability 1,
\begin{eqnarray}
\label{Ito-formula}
\lefteqn{f(Y(t))-f(Y(0))=\int_0^t f'(Y(s))G(s)ds+} \\
\nonumber
& & \int_0^t \int_{\bR^d}\int_{\{|z|>1\}} [f(Y(s-)+K(s,x,z))-f(Y(s-))]N(ds,dx,dz)+ \\
\nonumber
& & \int_0^t \int_{\bR^d}\int_{\{|z|\leq 1\}} [f(Y(s-)+H(s,x,z))-f(Y(s-))]\widehat{N}(ds,dx,dz)+\\
\nonumber
& & \int_0^t \int_{\bR^d}\int_{\{|z|\leq 1\}} [f(Y(s)+H(s,x,z))-f(Y(s))-H(s,x,z)f'(Y(s))]\nu(dz)dxds.
\end{eqnarray}
\end{theorem}

Note that since the first two terms on the right-hand side of \eqref{integral-proc} are processes of finite variation and the last term is a square-integrable martingale, $Y$ is a semimartingale. Therefore, the It\^o formula given by Theorem \ref{Ito-formula-th} can be derived from the corresponding result for a general semimartingale, assuming that $Y$ has c\`adl\`ag trajectories (see e.g. Theorem 2.5 of \cite{kunita04}).

The goal of the present article is to give an alternative proof of this result which contains the explicit construction of the c\`adl\`ag modification of $Y$ for which the It\^o formula holds.

We will also give the proof of the following variant of the It\^o formula, which will be useful for the applications related to the (finite-variance) L\'evy white noise, discussed in Section \ref{applic-section}.

\begin{theorem}[It\^o Formula II]
\label{Ito-formula-th2}
Let $Y=\{Y(t)\}_{t \geq 0}$ be a process defined by
\begin{equation}
\label{integral-proc-R0}
Y(t)=\int_0^t G(s)ds+\int_0^t \int_{\bR^d}\int_{\bR_0}H(s,x,z)\widehat{N}(ds,dx,dz), \quad t \geq 0,
\end{equation}
where $G$ and $H$ are predictable processes which satisfy \eqref{cond-G}, respectively \eqref{cond-H-R0}. Then there exists a c\`adl\`ag modification of $Y$ (denoted also by $Y$) such that for any $t>0$, with probability 1,
\begin{eqnarray*}
\lefteqn{f(Y(t))-f(Y(0))= \int_0^t f'(Y(s))G(s)ds+ }\\
& & \int_0^t \int_{\bR^d}\int_{\bR_0}[f(Y(s-)+H(s,x,z))-f(Y(s-))]\widehat{N}(ds,dx,dz)+\\
& & \int_0^t \int_{\bR^d}\int_{\bR_0}[f(Y(s)+H(s,x,z))-f(Y(s))-H(s,x,z)f'(Y(s))]
\nu(dz)dxds.
\end{eqnarray*}
\end{theorem}

The method that we use for proving Theorems \ref{Ito-formula-th} and \ref{Ito-formula-th2} is similar to the
one described in Section 4.4.2 of \cite{applebaum09} in the case of classical L\'evy processes, the difference being that in our case, $N$ is a PRM on $\bR_{+} \times \bR^d \times \bR_0$ instead of $\bR_{+} \times \bR_0$.
This method relies on a double ``interlacing'' technique, which
consists in first approximating the set $\{|z|\leq 1\}$ of small jumps by sets of the form $\{\varepsilon_n<|z| \leq 1\}$ with $\varepsilon_n \downarrow 0$ (in the case when $H$ and $K$ vanish outside a bounded Borel set $B \subset \bR^d$), and then approximating the spatial domain $\bR^d$ by regions of the form $[-a_n,a_n]^d$ with $a_n \uparrow \infty$. This approximation method is described in Section \ref{approx-section}.
Section \ref{proof-Ito-section} is dedicated to the proofs of Theorems \ref{Ito-formula-th} and \ref{Ito-formula-th2}. Finally, in Section \ref{applic-section} we discuss two applications of Theorem  \ref{Ito-formula-th2} in the case of the (finite-variance) L\'evy white noise introduced in \cite{B15}.

\section{Approximation by c\`adl\`ag processes}
\label{approx-section}

In this section, we show that the L\'evy-type integral processes given by \eqref{integral-proc} and \eqref{integral-proc-R0} have c\`adl\`ag modifications which are constructed by approximation. These modifications will play an important role in the proof of It\^o's formula. Since the process $Y_c(t)=\int_0^tG(s)ds$ is continuous, we assume that $G=0$.

We consider first processes of the form \eqref{integral-proc}.
We start by examining the case when both integrands $H$ and $K$ vanish outside a set $B \in \cB_b(\bR^d)$. Since the process
$\{\int_0^t \int_B \int_{\{|z|>1\}}K(s,x,z)N(ds,dx,dz);t \geq 0\}$
is clearly c\`adl\`ag (the integral being a sum with finitely many terms), we need to consider only the integral process which depends on $H$.

Note that if $H$ vanishes a.e. on $\Omega \times [0,T] \times B \times \{z \in \bR_0;|z|\leq \varepsilon\}$ for some $T>0$ and $\varepsilon\in (0,1)$, then
\begin{eqnarray*}
\lefteqn{\int_0^t \int_{B}\int_{\{|z|\leq 1\}}H(s,x,z)\widehat{N}(ds,dx,dz)=}\\
& & \int_0^t \int_{B}\int_{\{\varepsilon<|z|\leq 1\}}H(s,x,z)N(ds,dx,dz)-\int_0^t \int_{B}\int_{\{\varepsilon<|z|\leq 1\}}H(s,x,z)\nu(dz)dxds
\end{eqnarray*}
is a c\`adl\`ag process (the first term is a sum with finitely many terms and the second term in continuous).
Therefore, we will suppose that $H$ satisfies the following assumption:

\vspace{3mm}

{\em Assumption $A$.} It is not possible to find $T>0$ and $\varepsilon \in (0,1)$ such that
$$H(\omega,s,x,z)=0 \quad {\rm a.e. \ on} \quad \Omega \times [0,T] \times B \times \{z \in \bR_0;|z|\leq \varepsilon\}$$
with respect to the measure $P \times \mu$.

\begin{lemma}
\label{approx-th1}
Let $Y=\{Y(t)\}_{t \geq 0}$ be a process defined by
$$Y(t)=\int_0^t \int_{B}\int_{\{|z|\leq 1\}}H(s,x,z)\widehat{N}(ds,dx,dz),$$
where $B \in \cB_b(\bR^d)$ and $H$ is a predictable process which satisfies Assumption $A$ and
\begin{equation}
\label{cond-H-B}
E\int_0^t \int_{B}\int_{\{|z| \leq 1\}}|H(s,x,z)|^2 \nu(dz)dxds<\infty \quad \mbox{for all} \ t>0.
\end{equation}
Then, there exists a c\`adl\`ag modification $\widetilde{Y}=\{\widetilde{Y}(t)\}_{t \geq 0}$ of $Y$ such that for all $T>0$,
$$\sup_{t \leq T}|Y_n(t)-\widetilde{Y}(t)| \to 0 \quad a.s.,$$
where $$Y_n(t)=\int_0^t \int_B \int_{\{\varepsilon_n<|z| \leq 1\}}H(s,x,z) \widehat{N}(ds,dx,dz)$$ for some sequence $(\varepsilon_n)_n$ (depending on $T$) such that $\varepsilon_n \downarrow 0$ .
\end{lemma}

\noindent {\bf Proof:}  We use the same argument as in the proof of Theorem 4.3.4 of \cite{applebaum09}. Fix $T>0$. Let $\varepsilon_n=\sup\{\varepsilon>0;I(\varepsilon) \leq 8^{-n}\}$ where
$$I(\varepsilon)=E\int_0^T \int_{B}\int_{\{|z| \leq \varepsilon\}}|H(s,x,z)|^2\nu(dz)dxds.$$ Note that $(\varepsilon_{n})_n$ is non-increasing and $\varepsilon_n \downarrow 0$. (If $\varepsilon_n \downarrow \varepsilon_{*}>0$ then $I(\varepsilon_*) \leq I(\varepsilon_n) \leq 8^{-n}$ for all $n$. Hence $I(\varepsilon_*)=0$, which contradicts Assumption $A$.)

Note that $Y_n$ is a c\`adl\`ag martingale. By Doob's submartingale inequality and relation \eqref{E-int-hatN},
\begin{eqnarray*}
\lefteqn{E(\sup_{t \leq T}|Y_{n+1}(t)-Y_n(t)|^2)  \leq  4 E|Y_{n+1}(T)-Y_n(T)|^2 =} \\
& & 4E \int_0^T \int_{B} \int_{\{\varepsilon_{n+1} < |z| \leq \varepsilon_n\}} |H(s,x,z)|^2 \nu(dz)dxds \leq 4I(\varepsilon_n) \leq \frac{4}{8^n}.
\end{eqnarray*}
By Chebyshev's inequality,
$P(\sup_{t \leq T}|Y_{n+1}(t)-Y_n(t)|>2^{-n}) \leq 2^{-n+2}$.
By Borel-Cantelli lemma, with probability 1, the sequence $(Y_n)_n$ is Cauchy in the space $D[0,T]$ of c\`adl\`ag functions on $[0,T]$ equipped with the sup-norm. Its limit $\widetilde{Y}$ is a modification of $Y$ since for any $t \in [0,T]$, $\{Y_n(t)\}_n$ also converges to $Y(t)$ in $L^2(\Omega)$. Finally, we note that the process $\widetilde{Y}$ does not depend on $T$ (although the approximation sequence $(Y_n)_n$ does). If $\widetilde{Y}^{(T)}$ is the modification of $Y$ on $[0,T]$ and
$\widetilde{Y}^{(T')}$ is the modification of $Y$ on $[0,T']$ with $T<T'$, then $\widetilde{Y}^{(T)}(t)=\widetilde{Y}^{(T')}(t)$ a.s. for any $t \in [0,T]$. Hence, $\widetilde{Y}$ can be extended to $[0,\infty)$.
$\Box$

\vspace{3mm}

We consider now the case when the at least one of the integrands $H$ and $K$ do not vanish outside a set $B \in \cB_b(\bR^d)$. More precisely, we introduce the following assumptions:

{\em Assumption $B$.} It is not possible to find $T>0$ and $B \in \cB_b(\bR^d)$ such that $$H(\omega,t,x,z)=0 \quad {\rm a.e. \ on} \ \Omega \times [0,T] \times B^c \times \{z \in \bR_0;|z| \leq 1\}$$
with respect to the measure $P\times \mu$.

\vspace{2mm}

{\em Assumption $B'$.} It is not possible to find $T>0$ and $B \in \cB_b(\bR^d)$ such that $$K(\omega,t,x,z)=0 \quad {\rm a.e. \ on} \ \Omega \times [0,T] \times B^c \times \{z \in \bR_0;|z| > 1\}$$
with respect to the measure $P\times \mu$.

We consider bounded Borel sets in $\bR^d$ of the form $K_{a}=[-a,a]^d, a>0$.

\begin{theorem}[Interlacing I]
\label{approx-th2}
Let $Y=\{Y(t)\}_{t \geq 0}$ be a process defined by \eqref{integral-proc}
with $G=0$, where $H$ and $K$ are predictable processes which satisfy conditions \eqref{cond-H}, respectively \eqref{cond-K}, such that either $H$ satisfies Assumption $B$, or $K$ satisfies Assumption $B'$. Then, there exists a c\`adl\`ag modification $\widetilde{Y}=\{\widetilde{Y}(t)\}_{t \geq 0}$ of $Y$ such that for all $T>0$,
\begin{equation}
\label{Yn-approx-Y}
\sup_{t \leq T}|\widetilde{Y}_n(t)-\widetilde{Y}(t)| \to 0 \quad a.s.,
\end{equation}
where $\widetilde{Y}_n$ is a c\`adl\`ag modification of the process $Y_n$ defined by
$$Y_n(t)=\int_0^t \int_{E_n}\int_{\{|z|\leq 1\}}H(s,x,z)\widehat{N}(ds,dx,dz)+\int_0^t \int_{E_n}\int_{\{|z|>1\}}K(s,x,z)N(ds,dx,dz)$$
with $E_n=K_{a_n}$ for some sequence $(a_n)_n$ (depending on $T$) such that $a_n \uparrow \infty$.
\end{theorem}

\noindent {\bf Proof:} Fix $T>0$. Let $a_n=\inf\{a>0;I(a) \leq 8^{-n}\}$
where
$$I(a)=E\int_0^T\int_{K_a^c}\int_{\{|z|\leq 1\}}|H(s,x,z)|^2 \nu(dz)dxds+
E\int_0^T \int_{K_a^c}\int_{\{|z|>1\}}|K(s,x,z)|\nu(dz)dxds.$$
Note that $(a_{n})_n$ is non-decreasing and $a_n \uparrow \infty$. (If $a_n \uparrow a^{*}<\infty$ then $I(a^*) \leq I(a_n) \leq 8^{-n}$ for all $n$, and hence $I(a^*)=0$, which contradicts Assumptions $B$ or $B'$.)
Let $Y_n$ be the process given in the statement of the theorem with $E_n=K_{a_n}$. We denote by $Y_n^{(1)}(t)$ and $Y_n^{(2)}(t)$ the two integrals which compose $Y_n(t)$, depending on $H$, respectively $K$.

We denote by $\widetilde{Y}_n^{(1)}$ the c\`adl\`ag modification of $Y_n^{(1)}$ given by Lemma \ref{approx-th1}. By Doob's submartingale inequality and relation \eqref{E-int-hatN},
\begin{eqnarray*}
E(\sup_{t \leq T}|\widetilde{Y}_{n+1}^{(1)}(t)-\widetilde{Y}_n^{(1)}(t)|^2) & \leq & 4E\int_0^T \int_{E_{n+1} \verb2\2 E_n}\int_{\{|z|\leq 1\}}|H(s,x,z)|^2 \nu(dz)dxds \\
& \leq & 4 I(a_n) \leq \frac{4}{8^n}.
\end{eqnarray*}
By Chebyshev's inequality,
$P(\sup_{t \leq T}|\widetilde{Y}_{n+1}^{(1)}(t)-\widetilde{Y}_n^{(1)}(t)|>
2^{-n-1}) \leq 2^{-n+4}$.

Note that $Y_n^{(2)}$ is a c\`adl\`ag process. For any $t \in [0,T]$,
$$|Y_{n+1}^{(2)}(t)-Y_{n}^{(2)}(t)| \leq \int_0^t \int_{E_{n+1} \verb2\2 E_n}\int_{\{|z|>1\}}|K(s,x,z)|N(ds,dx,dz),$$
and hence, using relation \eqref{E-int-N},
\begin{eqnarray*}
E(\sup_{t \leq T}|Y_{n+1}^{(2)}(t)-Y_{n}^{(2)}(t)|) & \leq & E\int_0^T \int_{E_{n+1} \verb2\2 E_n}\int_{\{|z|>1\}}|K(s,x,z)|\nu(dz)dxds \\
& \leq & I(a_n) \leq \frac{1}{8^n}.
\end{eqnarray*}
By Markov's inequality,
$P(\sup_{t \leq T}|Y_{n+1}^{(2)}(t)-Y_{n}^{(2)}(t)|>2^{-n-1}) \leq 2^{-2n+1}$.

Let $\widetilde{Y}_n(t)=\widetilde{Y}_n^{(1)}(t)+Y_n^{(2)}(t)$. Then
$P(\sup_{t \leq T}|\widetilde{Y}_{n+1}(t)-\widetilde{Y}_{n}(t)|>2^{-n}) \leq 2^{-n+4}+2^{-2n+1}$, and the conclusion follows by the Borel-Cantelli Lemma, as in the proof of Lemma \ref{approx-th1}. $\Box$

\vspace{3mm}

We consider next processes of the form \eqref{integral-proc-R0} with $G=0$. Note that if $H$ vanishes a.e. outside a set $B \in \cB_b(\bR^d)$ then
\begin{eqnarray*}
Y(t)&=&\int_0^t \int_{B}\int_{\{|z|\leq 1\}}H(s,x,z)\widehat{N}(ds,dx,dz)+\int_0^t \int_{B}\int_{\{|z|>1\}}H(s,x,z)N(ds,dx,dz)\\
& & -\int_0^t \int_{B}\int_{\{|z|>1\}}H(s,x,z)\nu(dz)dxds,
\end{eqnarray*}
where the first term has a c\`adl\`ag modification given by Lemma \ref{approx-th1}, the second term is c\`adl\`ag, and the third term is continuous.
Therefore, we will suppose that $H$ satisfies the following assumption:

\vspace{3mm}

{\em Assumption $C$.} It is not possible to find $T>0$ and $B \in \cB_b(\bR^d)$ such that
$$H(\omega,s,x,z)=0 \quad {\rm a.e.} \quad \Omega \times [0,T] \times B^c \times \bR_0$$
with respect to the measure $P \times \mu$.

\begin{theorem}[Interlacing II]
\label{approx-th-R0}
Let $Y$ be a process given by \eqref{integral-proc-R0} with $G=0$, where $H$ is a predictable process which satisfies \eqref{cond-H-R0} and Assumption $C$.
Then, there exists a c\`adl\`ag modification $\widetilde{Y}=\{\widetilde{Y}(t)\}_{t \geq 0}$ of $Y$ such that
\eqref{Yn-approx-Y} holds, where $\widetilde{Y}_{n}$ is a c\`adl\`ag modification of the process $Y_n$ defined by:
$$Y_n(t)=\int_0^t \int_{E_n}\int_{\bR_0}H(s,x,z)\widehat{N}(ds,dx,dz),$$
with $E_n=K_{a_n}$ for some sequence $(a_n)_n$ (depending on $T$) such that $a_n \uparrow \infty$.
\end{theorem}

\noindent {\bf Proof:} We proceed as in the proof of Theorem \ref{approx-th2}. Fix $T>0$. Let $a_n=\inf\{a>0;I(a) \leq 8^{-n}\}$ where
$$I(a)=\int_0^t \int_{K_a^c}\int_{\bR_0}|H(s,x,z)|^2 \nu(dz)dxds.$$
By Assumption $C$, $a_n \uparrow \infty$. We write $Y_n(t)$ as the sum of two integrals, corresponding to the regions $\{|z|\leq 1\}$, and $\{|z|>1\}$. We denote these integrals by $Y_n^{(1)}(t)$, respectively $Y_n^{(2)}(t)$. Note that $Y_n^{(2)}$ is c\`adl\`ag. Let $\widetilde{Y}_n^{(1)}$ be the c\`adl\`ag modification of $Y_n^{(1)}$ given by Lemma \ref{approx-th1}.

Let $\widetilde{Y}_n(t)=\widetilde{Y}_n^{(1)}(t)+Y_n^{(2)}(t)$. By Doob's submartingale inequality,
$$E(\sup_{t \leq T}|\widetilde{Y}_{n+1}(t)-\widetilde{Y}_n(t)|^2)\leq 4
E\int_0^T\int_{E_{n+1}\verb2\2 E_n}\int_{\bR_0}|H(s,x,z)|^2 \nu(dz)dxds$$
and the conclusion follows as in the proof of Lemma \ref{approx-th1}. $\Box$

\section{Proof of It\^o Formula}
\label{proof-Ito-section}

In this section, we give the proofs of Theorem \ref{Ito-formula-th} and Theorem \ref{Ito-formula-th2}.

We start with the simpler case when there are no small jumps (the analogue of Lemma 4.4.6 of \cite{applebaum09}).

\begin{lemma}
\label{Ito-formula1}
Let
$$Y(t)=\int_0^tG(s)ds+\int_{0}^{t}\int_{B}\int_{\{|z|>\varepsilon\}}
K(s,x,z)N(ds,dx,dz)=:Y_c(t)+Y_d(t),$$
where $G$ is a predictable process which satisfies \eqref{cond-G}, $B \in \cB_b(\bR^d)$, $\varepsilon>0$ and $K$ is a predictable process. Then, for any function $f\in C^1(\bR)$ and for any $t>0$,
\begin{eqnarray*}
\lefteqn{f(Y(t))-f(Y(0))=\int_0^t f'(Y(s))G(s)ds+}\\
& & \int_0^t \int_{B}\int_{\{|z|>\varepsilon\}}[f(Y(s-)+K(s,x,z))-f(Y(s-))]N(ds,dx,dz).
\end{eqnarray*}
\end{lemma}

\noindent {\bf Proof:} We denote $\Gamma=\{|z|>\varepsilon\}$. By Proposition 5.3 of \cite{resnick07}, we may assume that the restriction of $N$ to the set $\bR_{+} \times B \times \Gamma$ has points $(T_i,X_i,Z_i),i \geq 1$, where $T_1<T_2< \ldots$ are the points of a Poisson process on $\bR_{+}$ of intensity $\lambda=|B|\nu(\Gamma)$ and $\{(X_i,Z_i)\}_{i \geq 1}$ are i.i.d. on $B \times \Gamma$ with distribution $\lambda^{-1}dx\nu(dz)$, independent of $(T_i)_{i \geq 1}$. We consider two cases.

{\em Case 1: $G=0$.} By the representation of $N$,
$Y(t)=\sum_{T_i \leq t}K(T_i,X_i,Z_i)$.
So $t \mapsto Y(t)$ is a step function which has a jump of size $K(T_i,X_i,Z_i)$ at each point $T_i$ and $Y(T_i-)=Y(T_{i-1})$. Hence
\begin{eqnarray*}
f(Y(t))-f(Y(0))&=& \sum_{T_i \leq t}[f(Y(T_i))-f(Y(T_{i-1}))]\\
&=& \sum_{T_i \leq t}[f(Y(T_{i}-)+K(T_i,X_i,Z_i))-f(Y(T_{i}-))],
\end{eqnarray*}
and the conclusion follows since $N$ has points $(T_i,X_i,Z_i)$ in $\bR_{+} \times B \times \Gamma$.

{\em Case 2: $G$ is arbitrary.} The map $t \mapsto Y_d(t)$ is a step function which has a jump of size $K(T_i,X_i,Z_i)$ at time $T_i$. Since $Y_c$ is continuous, the jump times and the jump sizes of $Y$ coincide with those of $Y_d$, i.e. $\Delta Y(T_i)=\Delta Y_d(T_i)=K(T_i,X_i,Z_i)$.
We use the decomposition $$f(Y(t))-f(Y(0))=A(t)+B(t),$$
where $A$ and $B$ are defined as follows: if $T_{n-1} \leq t<T_n$, we let
\begin{eqnarray*}
A(t)&=&\sum_{i=1}^{n-1}[f(Y(T_i))-f(Y(T_{i}-))] \\
B(t)&=&\sum_{i=1}^{n-1}[f(Y(T_{i}-))-f(Y(T_{i-1}))]+[f(Y(t))-f(Y(T_{n-1}))].
\end{eqnarray*}

\noindent Note that
\begin{eqnarray*}
A(t)&=& \sum_{i=1}^{n}[f(Y(T_{i}-)+K(T_i,X_i,Z_i))-f(Y(T_{i}-))] \\
&=& \int_0^t \int_{B}\int_{\Gamma}[f(Y(s-)+K(s,x,z))-f(Y(s-))]N(ds,dx,dz).
\end{eqnarray*}

It remains to prove that
\begin{equation}
\label{equality-B}
B(t)=\int_0^t f'(Y(s))G(s)ds.
\end{equation}
For this, we assume that $T_{n-1} \leq t <T_n$ and we write
$$\int_0^t f'(Y(s))G(s)ds=\sum_{i=1}^{n-1}\int_{T_{i-1}}^{T_i}f'(Y(s))G(s)ds+\int_{T_{n-1}}^{t}
f'(Y(s))G(s)ds.$$
So it suffices to prove that
\begin{equation}
\label{equality-B1}
\int_{T_{i-1}}^{T_i}f'(Y(s))G(s)ds=f(Y(T_{i}-))-f(Y(T_{i-1}))
\end{equation}
for all $i=1,\ldots,n-1$, and
\begin{equation}
\label{equality-B2}
\int_{T_{n-1}}^{t}f'(Y(s))G(s)ds=f(Y(t))-f(Y(T_{n-1})).
\end{equation}

We first prove \eqref{equality-B1}. Fix $i=1,\ldots,n-1$. For any $s \in (T_{i-1},T_i)$, $Y(s)=Y_c(s)+Y_d(T_{i-1}):=g_i(s)$ and $g_i'(s)=Y_c'(s)=G(s)$. We extend $g_i$ by continuity to $[T_{i-1},T_i]$. Hence
\begin{eqnarray*}
\int_{T_{i-1}}^{T_i}f'(Y(s))G(s)ds &=& \int_{T_{i-1}}^{T_i}
f'(g_i(s))g_i'(s)ds=f(g_i(T_i))-f(g_i(T_{i-1}))\\
&=& f(Y_c(T_i)+Y_d(T_{i-1}))-f(Y_c(T_{i-1})+Y_d(T_{i-1})) \\
&=&  f(Y(T_i-))-f(Y(T_{i-1})),
\end{eqnarray*}
where for the last equality we used the fact that $Y_d(T_{i-1})=Y_d(T_{i}-)$ and hence $Y_c(T_i)+Y_d(T_{i-1})=Y_c(T_i-)+Y_d(T_i-)=Y(T_i-)$. This proves \eqref{equality-B1}.

Next, we prove \eqref{equality-B2}. Note that if $t=T_{n-1}$, both terms are zero. So, we assume that $t>T_{n-1}$. For any $s \in (T_{n-1},t)$, $Y(s)=Y_c(s)+Y_{d}(T_{n-1}):=g(s)$ and $g'(s)=Y_c'(s)=G(s)$. Arguing as above, we see that
\begin{eqnarray*}
\int_{T_{n-1}}^{t}f'(Y(s))G(s)ds &=& \int_{T_{n-1}}^{t}
f'(g(s))g'(s)ds=f(g(t))-f(g(T_{n-1}))\\
&=& f(Y_c(t)+Y_d(T_{n-1}))-f(Y_c(T_{n-1})+Y_d(T_{n-1})) \\
&=&  f(Y(t))-f(Y(T_{n-1})),
\end{eqnarray*}
where for the last equality we used the fact that $Y_d(T_{n-1})=Y_d(t)$ and hence $Y_c(t)+Y_d(T_{n-1})=Y_c(t)+Y_d(t)=Y(t)$. This concludes the proof of  \eqref{equality-B2}.
$\Box$

\vspace{3mm}

\noindent {\bf Proof of Theorem \ref{Ito-formula-th}:}  We fix $t>0$. We assume that $f'$ and $f''$ are bounded. (Otherwise, we use $\tau_k=\inf\{s>0; |Y(s)|>k\}$ for $k \geq 1$.)

 {\em Case 1: $H$ and $K$ vanish outside a fixed set $B \in \cB_b(\bR^d)$}.

If $H$ vanishes a.e. on $\Omega \times [0,T] \times B \times \{z \in \bR_0;|z| \leq \varepsilon\}$ for some $T>0$ and $\varepsilon \in (0,1)$, the conclusion follows from Lemma \ref{Ito-formula1}. Therefore, we suppose that $H$ satisfies Assumption $A$. By Lemma \ref{approx-th1}, there exists a c\`adl\`ag modification of $Y$ (denoted also by $Y$) such that
\begin{equation}
\label{sup-converges}
\sup_{s \leq t}|Y_n(s)-Y(s)| \to 0,
\end{equation}
where the process $\{Y_n(s)\}_{s \in [0,t]}$ is defined by
\begin{eqnarray*}
Y_n(s)&=&\int_0^s G(r)dr+\int_0^s \int_{B}\int_{\{\varepsilon_n<|z|\leq 1\}}H(r,x,z)\widehat{N}(dr,dx,dz)+\\
& & \int_0^s \int_{B}\int_{\{|z|>1\}}K(r,x,z)N(dr,dx,dz), \quad s \in [0,t],
\end{eqnarray*}
$(\varepsilon_n)_n$ being the sequence given by Lemma \ref{approx-th1} with $T=t$. Consequently,
\begin{equation}
\label{sup-converges-ll}
\sup_{s \leq t}|Y_n(s-)-Y(s-)| \to 0.
\end{equation}
 Note that
$$Y_n(s)=\int_0^s \overline{G}(r)dr+\int_0^s \int_{B}\int_{\{|z|>\varepsilon_n\}}\overline{K}(r,x,z)N(rs,dx,dz),$$
where $\overline{G}(s)=G(s)-\int_{B}\int_{\{\varepsilon_n<|z| \leq 1\}}H(s,x,z)\nu(dz)dx$ and $\overline{K}(s,x,z)=H(s,x,z)\linebreak 1_{\{|z|\leq 1\}}+K(s,x,z)
1_{\{|z|>1\}}$. By the Cauchy-Schwarz inequality, $\overline{G}$ satisfies \eqref{cond-G} (since $B$ is a {\em bounded} set and $H$ satisfies \eqref{cond-H-B}).
We apply Lemma \ref{Ito-formula1} to $Y_n$:
\begin{eqnarray*}
\lefteqn{f(Y_n(t))-f(Y_n(0))=\int_0^t f'(Y_n(s))\overline{G}(s)ds+}\\
& & \int_0^t \int_{B}\int_{\{|z|>\varepsilon_n\}}[f(Y_n(s-)+\overline{K}(s,x,z))-f(Y_n(s-))]
N(ds,dx,dz).
\end{eqnarray*}
After using the definitions of $\overline{G}$ and $\overline{K}$, as well as adding and subtracting $$\int_0^t \int_{B}\int_{\{\varepsilon_n<|z|\leq 1\}}[f(Y_n(s)+H(s,x,z))-f(Y_n(s))]\nu(dz)dxds,$$ we obtain that:
\begin{eqnarray}
\nonumber
\lefteqn{f(Y_n(t))-f(Y_n(0))=  \int_0^t f'(Y_n(s))G(s)ds+}\\
\nonumber
& & \int_0^t \int_{B}\int_{\{|z|>1\}} [f(Y_n(s-)+K(s,x,z))-f(Y_n(s-))]N(ds,dx,dz)+ \\
\nonumber
& & \int_0^t \int_{B}\int_{\{|z|\leq 1\}} [f(Y_n(s-)+H(s,x,z))-f(Y_n(s-))]\widehat{N}(ds,dx,dz)+\\
\nonumber
& &  \int_0^t \int_{B}\int_{\{|z|\leq 1\}} [f(Y_n(s)+H(s,x,z))-f(Y_n(s))-H(s,x,z)\\
\label{4terms}
& & \ \ \ \ \ \ \ \ \ \ \ \ \ \ \ \ \ f'(Y_n(s))]\nu(dz)dxds:=T_{1,n}+T_{2,n}+T_{3,n}+T_{4,n}.
\end{eqnarray}

We denote by $T_1,T_2,T_3$, respectively $T_4$ the four terms on the right-hand side of \eqref{Ito-formula}.
The conclusion will follow by taking the limit as $n \to \infty$ in \eqref{4terms}.
The left-hand side converges to $f(Y(t))-f(Y(0))$, by \eqref{sup-converges}.

We treat separately the four terms in the right-hand side.
By the dominated convergence theorem,
$$E|T_{1,n}-T_1| \leq E\int_0^t |f'(Y_n(s))-f'(Y(s))||G(s)|ds \to 0.$$
Since $T_{2,n}$ is a sum with a finite number of terms, using \eqref{sup-converges} and the continuity of $f$, we see that $T_{2,n} \to T_2$ a.s. For the third term, note that $E|T_{3,n}-T_3|^2 \leq 2(A_n+B_n)$, where
\begin{eqnarray*}
A_n&=& E\int_0^t \int_B \int_{\{\varepsilon_n<|z| \leq 1\}}|V_n(s,x,z)-V(s,x,z)|^2 \nu(dz)dxds,\\
B_n&=& E\int_0^t \int_{B}\int_{\{|z| \leq \varepsilon_n\}}|V(s,x,z)|^2 \nu(dz)dxds,
\end{eqnarray*}
and $V_n(s,x,z):=f(Y_n(s)+H(s,x,z))-f(Y_n(s)) \to
V(s,x,z):=f(Y(s) +H(s,x,z))-f(Y(s))$ a.s., by \eqref{sup-converges} and the continuity of $f$. By the dominated convergence theorem, $A_n \to 0$ and $B_n \to 0$. To justify the application of this theorem, we use Taylor's formula of the first order:
\begin{equation}
\label{Taylor-1}
f(b)-f(a)=(b-a)\int_0^1 f'(a+\theta(b-a))d\theta,
\end{equation}
and the fact that $f'$ is bounded. This proves that $T_{3,n} \to T_3$ in $L^2(\Omega)$.

Finally, $E|T_{4,n}-T_4| \leq C_n+D_n$, where
\begin{eqnarray*}
C_n&=&E \int_0^t \int_{B} \int_{\{\varepsilon_n<|z|\leq 1\}}|U_n(s,x,z)-U(s,x,z)|\nu(dz)dxds, \\
D_n&=&E\int_0^t \int_B \int_{|z|\leq \varepsilon_n}|U(s,x,z)|\nu(dz)dxds,
\end{eqnarray*}
and $U_n(s,x,z):=f(Y_n(s)+H(s,x,z))-f(Y_n(s))-H(s,x,z)f'(Y_n(s))  \to
U(s,x,z):=f(Y(s)+H(s,x,z))-f(Y(s))-H(s,x,z)f'(Y_n(s))$ a.s., by \eqref{sup-converges-ll} and the continuity of $f$. By the dominated convergence theorem, $C_n \to 0$ and $D_n \to 0$. To justify the application of this theorem, we use Taylor's formula of second order:
\begin{equation}
\label{Taylor-2}
f(b)-f(a)=(b-a)f'(a)+(b-a)^2\int_0^1 f''(a+\theta(b-a))(1-\theta)d\theta,
\end{equation}
and the fact that $f''$ is bounded. This proves that $T_{4,n} \to T_4$ in $L^1(\Omega)$.

{\em Case 2. $H$ satisfies Assumption $B$ or $K$ satisfies Assumption $B'$.}\\
By Theorem \ref{approx-th2}, there exists a c\`adl\`ag approximation of $Y$ (denoted also by $Y$) such that
\eqref{sup-converges} holds, where
$\{Y_n(s)\}_{s \in [0,t]}$ is a c\`adl\`ag modification of
\begin{eqnarray*}
Y_n(s)&=&\int_0^s G(r)dr+\int_0^s \int_{E_n}\int_{\{|z| \leq 1\}}H(r,x,z)\widehat{N}(dr,dx,dz)+\\
& & \int_0^s \int_{E_n}\int_{\{|z|>1\}}K(r,x,z)N(dr,dx,dz), \quad s \in [0,t],
\end{eqnarray*}
$(E_n)_n \subset \cB_b(\bR^d)$ being the sequence given by Theorem \ref{approx-th2} with $T=t$. Using the result of Case 1 for the process $Y_n$, we obtain
\begin{eqnarray*}
\lefteqn{f(Y_n(t))-f(Y_n(0))=\int_0^t f'(Y_n(s))G(s)ds+}\\
& & \int_0^t \int_{E_n}\int_{\{|z|>1\}} [f(Y_n(s-)+K(s,x,z))-f(Y_n(s-))]N(ds,dx,dz) \\
& & + \int_0^t \int_{E_n}\int_{\{|z|\leq 1\}} [f(Y_n(s-)+H(s,x,z))-f(Y_n(s-))]\widehat{N}(ds,dx,dz)\\
& & +  \int_0^t \int_{E_n}\int_{\{|z|\leq 1\}} [f(Y_n(s)+H(s,x,z))-f(Y_n(s))-H(s,x,z)f'(Y_n(s)]\nu(dz)dxds.\\
\end{eqnarray*}
The conclusion follows letting $n \to \infty$ as in Case 1. $\Box$

\vspace{3mm}

\noindent {\bf Proof of Theorem \ref{Ito-formula-th2}:} We assume that $f'$ and $f''$ are bounded. We fix $t$.

{\em Case 1. $H$ vanishes outside a set $B \in \cB_b(\bR^d)$.}
We write
\begin{eqnarray*}
Y(t)&=& \int_0^t \overline{G}(s)ds+\int_0^t \int_{B} \int_{\{|z|\leq 1\}}H(s,x,z)\widehat{N}(ds,dx,dz)+\\
& & \int_0^t \int_{B} \int_{\{|z|>1\}}H(s,x,z)N(ds,dx,dz),
\end{eqnarray*}
where $\overline{G}(s)=G(s)-\int_{B}\int_{\{|z|>1\}}H(s,x,z)\nu(dz)dx$. By the Cauchy-Schwarz inequality, $\overline{G}$ satisfies \eqref{cond-G} (since $B$ is a {\em bounded} set). By Theorem \ref{Ito-formula-th}, there exists a c\`adl\`ag modification of $Y$ (denoted also by $Y$) such that
\begin{eqnarray*}
\lefteqn{f(Y(t))-f(Y(0))=\int_0^t f'(Y(s))\overline{G}(s)ds+ } \\
& & \int_0^t \int_B \int_{\{|z|>1\}}[f(Y(s-)+H(s,x,z))-f(Y(s-))]N(ds,dx,dz)+\\
& & \int_0^t \int_{B} \int_{\{|z|\leq 1\}} [f(Y(s-)+H(s,x,z))-f(Y(s-))]\widehat{N}(ds,dx,dz)+\\
& & \int_0^t \int_{B}\int_{\{|z|\leq 1\}}[f(Y(s)+H(s,x,z))-f(Y(s))-H(s,x,z)f'(Y(s))]\nu(dz)dxds.
\end{eqnarray*}
We add and subtract $\int_0^t \int_{B} \int_{\{|z|>1\}}[f(Y(s)+H(s,x,z))-f(Y(s))]\nu(dz)dxds$. The conclusion follows by rearranging the terms.

{\em Case 2. $H$ satisfies Assumption $C$.}\\
By Theorem \ref{approx-th-R0}, there exists a c\`adl\`ag modification of $Y$ (denoted also by $Y$) such that \eqref{sup-converges} holds, where $\{Y_n(s)\}_{s \in [0,t]}$ is a c\`adl\`ag modification of
$$Y_n(s)=\int_0^sG(r)dr+\int_0^s \int_{E_n}\int_{\bR_0}H(r,x,z)\widehat{N}(dr,dx,dz), \quad s \in [0,t],$$
$(E_n)_n$ being the sequence given by Theorem \ref{approx-th-R0} with $T=t$.
We write the It\^o formula for the process $Y_n$ (using Case 1) and we let $n \to \infty$.  $\Box$

\section{Applications}
\label{applic-section}

In this section, we assume that the L\'evy measure $\nu$ satisfies the condition:
$$v:=\int_{\bR_0}z^2 \nu(dz)<\infty.$$
 As in \cite{B15}, we consider the process $L=\{L(B); t \geq 0, B \in \cB_b(\bR_{+} \times \bR^d)\}$ defined by:
$$L(B)=\int_{B \times \bR_0}z\widehat{N}(ds,dx,dz).$$
For any predictable process $X=\{X(t,x);t \geq 0,x \in \bR^d\}$
such that
\begin{equation}
\label{cond-X}
E\int_0^T \int_{\bR^d}|X(t,x)|^2 dxdt<\infty \quad \mbox{for any} \ T>0,
\end{equation}
we can define the stochastic integral of $X$ with respect to $L$ and this integral satisfies:
$$\int_0^T \int_{\bR^d}X(t,x)L(dt,dx)=\int_0^T \int_{\bR^d}\int_{\bR_0}X(t,x)z\widehat{N}(dt,dx,dz).$$
By \eqref{E-int-hatN}, this integral has the following isometry property:
$$E\left|\int_0^T \int_{\bR^d}X(t,x)L(dt,dx)\right|^2=v E\int_{0}^{T}\int_{\bR^d}|X(t,x)|^2dxdt.$$

When used as a noise process perturbing an SPDE, $L$ behaves very similarly to the Gaussian white noise. For this reason, $L$ was called a {\em L\'evy white noise} in \cite{B15}.

\subsection{Kunita Inequality}

The following maximal inequality is due to Kunita (see Theorem 2.11 of \cite{kunita04}). In problems related to SPDEs with noise $L$, this result plays the same role as the Burkholder-Davis-Gundy inequality for SPDEs with Gaussian white noise.

\begin{theorem}[Kunita Inequality]
\label{Kunita-ineq}
Let $Y=\{Y(s)\}_{t \geq 0}$ be a process given by
$$Y(t)=\int_0^t\int_{\bR^d}X(s,x)L(ds,dx), \quad t \geq 0,$$ where $X$ is a predictable process which satisfies \eqref{cond-X}.

If
$m_p=\int_{\bR_0}|z|^p \nu(dz)<\infty$ for some $p \geq 2$, then for any $t>0$,
\begin{eqnarray*}
E(\sup_{s \leq t}|Y(s)|^p) \leq C_p \left\{E \left(\int_0^t \int_{\bR^d} |X(s,x)|^2 dxds \right)^{p/2}+ E\int_0^t \int_{\bR^d} |X(s,x)|^p dxds \right\},
\end{eqnarray*}
where $C_p=K_p \max(v^{p/2},m_p)$ and $K_p$ is the constant in Theorem 2.11 of \cite{kunita04}.
\end{theorem}

\noindent {\em Proof:} We apply Theorem \ref{Ito-formula-th2} with $f(x)=|x|^p$ and $H(s,x,z)=X(s,x)z$. The proof is identical to that of Theorem 2.11 of \cite{kunita04}. We omit the details. $\Box$

\begin{remark}
{\rm Kunita's constant $K_p$ cannot be computed explicitly. Theorem \ref{Kunita-ineq} is proved in \cite{BN15} using a different method which shows that $K_p$ is directly related to the constant $B_p$ in Rosenthal's inequality, which is $O(p/\ln p)$. }
\end{remark}

\subsection{It\^o Representation Theorem and Chaos Expansion}

In this section, we give an application to Theorem \ref{Ito-formula-th2} to exponential martingales, which leads to It\^o representation theorem and a chaos expansion (similarly to Sections 5.3-5.4 of \cite{applebaum09}).

For any $h \in L^2(\bR_{+} \times \bR^d)$ we let
 $L_h(t)=\int_0^t \int_{\bR^d}h(s,x)L(ds,dx)$ for $t\geq 0$. We work with the c\`adl\`ag modification of the process $L_h$ given by Theorem \ref{approx-th-R0}.
By Lemma 2.4 of \cite{B15},
$$E(e^{iL_h(t)})=\exp \left\{\int_0^t \int_{\bR^d}\Psi(h(s,x))dxds \right\},$$
where
$$\Psi(u)=\int_{\bR_0}(e^{iuz}-1-iuz)\nu(dz), \quad u \in \bR.$$
Hence $E(M_h(t))=1$ for all $t \geq 0$, where
$$M_h(t)=\exp \left\{iL_h(t)-\int_0^t \int_{\bR^d}\Psi(h(s,x))dxds\right\}, \quad t\geq 0.$$

The following result is the analogue of Lemma 5.3.3 of \cite{applebaum09}.

\begin{lemma}
\label{expo-mart-lemma}
For any $h \in L^2(\bR_{+} \times \bR^d)$ and $t>0$, with probability 1,
$$M_h(t)=1+\int_0^t \int_{\bR^d}\int_{\bR_0}(e^{ih(s,x)z}-1)M_h(s-)\widehat{N}(ds,dx,dz).$$
\end{lemma}

\noindent {\bf Proof:} We apply Theorem \ref{Ito-formula-th2} to the function $f(x)=e^{ix}$ and the process
$$Y(t)=L_h(t)+i\int_0^t \int_{\bR^d} \Psi(h(s,x))dxds.$$
Hence, $H(s,x,z)=h(s,x)z$ and $G(s)=i\int_{\bR^d}\Psi(h(s,x))dx$. We obtain:
\begin{eqnarray*}
\lefteqn{M_h(t)-1=\int_0^t \int_{\bR^d}\int_{\bR_0} (e^{iY(s-)+ih(s,x)z}-e^{iY(s-)})\widehat{N}(ds,dx,dz)+ }\\
& & \int_0^t \int_{\bR^d} \int_{\bR_0} (e^{iY(s)+ih(s,x)z}-e^{iY(s)}-izh(s,x)e^{iY(s)})\nu(dz)dxds+ \\
& & \int_0^t ie^{iY(s)} \left(i\int_{\bR^d}\Psi(h(s,x))dx\right)ds.
\end{eqnarray*}
Since the sum of the last two integrals is 0, the conclusion follows.  $\Box$

\vspace{3mm}

We fix $T>0$. We let $\cF_t^{L}=\sigma(\{L_s(B);0 \leq s \leq t,B \in \cB_b(\bR^d)\})$. We denote by $L_{\bC}^2(\Omega,\cF_T^L,P)$ be the space of $\bC$-valued square-integrable random variables which are measurable with respect to $\cF_T^L$.

\begin{lemma}
\label{dense-lemma}
The linear span of the set $\cA=\{M_h(T); h \in L^2(\bR_{+} \times \bR^d)\}$ is dense in $L_{\bC}^2(\Omega,\cF_T^L,P)$.
\end{lemma}

\noindent {\bf Proof:} The proof is similar to that of Lemma 5.3.4 of \cite{applebaum09}. We omit the details. $\Box$

\vspace{3mm}

%We have the following result.

\begin{theorem}[It\^o Representation Theorem]
\label{Ito-rep-th}
For any $F \in L_{\bC}^2(\Omega,\cF_T^L,P)$, there exists a unique predictable $\bC$-valued process $\psi=\{\psi(t,x,z); t \in [0,T],\linebreak x \in \bR^d,z\in \bR_0\}$ satisfying
\begin{equation}
\label{cond-psi}
E\int_0^T \int_{\bR^d}\int_{\bR_0}|\psi(t,x,z)|^2 \nu(dz)dxdt<\infty
\end{equation}
such that
\begin{equation}
\label{rep}
F=E(F)+\int_0^T \int_{\bR^d}\int_{\bR_0}\psi(t,x,z)\widehat{N}(dt,dx,dz).
\end{equation}
\end{theorem}

\noindent {\bf Proof:} By Lemma \ref{expo-mart-lemma}, relation \eqref{rep} holds for $F=M_h(T)$ with $\psi(t,x,z)=(e^{ih(t,x)z}-1)M_h(t-)$. The conclusion follows by an approximation argument using Lemma \ref{dense-lemma}. $\Box$

\vspace{3mm}

The multiple (and iterated) integral with respect $\widehat{N}$ can be defined similarly to the Gaussian white-noise case (see e.g. Section 5.4 of \cite{applebaum09}).

More precisely, we consider the Hilbert space $\cH=L^2(S, \cS,\mu)$, where $S=[0,T] \times \bR^d \times \bR_0$, $\cS=\cB([0,T]) \times \cB(\bR^d) \times \cB(\bR_0)$ and $\mu=dtdx \nu(dz)$. For any integer $n \geq 1$, we consider the $n$-th tensor product space $\cH^{\otimes n}=L^2(S^n,\cS^n, \mu^n)$.
The $n$-th multiple integral $I_n(f)$ with respect to $\widehat{N}$ can be constructed
for any function $f \in \cH^{\otimes n}$, and this integral has the isometry property:
$$E|I_n(f)|^2=n!\|f\|_{\cH^{\otimes n}}^{2}.$$
Moreover, if $n \not=m$, then $E[I_n(f)I_m(g)]=0$ for all $f\in \cH^{\otimes n}$ and $g\in \cH^{\otimes m}$.

\vspace{3mm}

It\^o representation theorem leads to the following result.

\begin{theorem}[Chaos Expansion]
For any $F \in L^2(\Omega,\cF_T^L,P)$, there exist some symmetric functions $f_n \in \cH^{\otimes n}$, $n \geq 1$ such that
$$F=E(F)+\sum_{n \geq 1}I_n(f_n) \quad \mbox{in} \ L^2(\Omega).$$
In particular,
$$E|F|^2=|E(F)|^2+\sum_{n \geq 1}n! \|f_n\|_{\cH^{\otimes n}}^{2}.$$
\end{theorem}

\noindent {\bf Proof:} We use the same argument as in the classical case, when $\widehat{N}$ is a PRM on $\bR_{+} \times \bR_0$ and $L(t)=\int_0^t \int_{\bR_0}z \widehat{N}(ds,dz), t \geq 0$ is a square-integrable L\'evy process (see Theorem 5.4.6 of \cite{applebaum09} or Theorem 10.2 of \cite{DOP09}). %To simplify the notation, we denote $u=(t,x,z)$.
By Theorem \ref{Ito-rep-th}, there exists a predictable process $\psi_1$ satisfying \eqref{cond-H-R0} such that
\begin{equation}
\label{rep-F-psi1}
F=E(F)+\int_0^T \int_{\bR^d}\int_{\bR_0}\psi_1(t_1,x_1,z_1)\widehat{N}(dt_1,dx_1,dz_1).
\end{equation}
By \eqref{cond-psi}, $E|\psi_1(t_1,x_1,z_1)|^2<\infty$ for almost all $(t_1,x_1,z_1)$.
For such $(t_1,x_1,z_1)$ fixed, we apply Theorem \ref{Ito-rep-th} again to the variable $\psi_1(t_1,x_1,z_1)$. Hence, there exists a predictable process $\psi_2=\{\psi_2(t_2,x_2,z_2);t_2 \in [0,t_1],x_2 \in \bR^d,z_2 \in \bR_0\}$ satisfying $E\int_0^{t_1} \int_{\bR^d}\int_{\bR_0}|\psi_1(t_1,x_1,z_1)|^2\nu(dz_1)dx_1dt_1<\infty$
such that
$$\psi_1(t_1,x_1,z_1)=E(\psi_1(t_1,x_1,z_1))+\int_0^{t_1} \int_{\bR^d}\int_{\bR_0}\psi_2(t_2,x_2,z_2)\widehat{N}(dt_2,dx_2,dz_2).$$
We substitute this into \eqref{rep-F-psi1} and iterate the procedure. We omit the details. $\Box$

\end{document}